\documentclass[reqno]{amsart}
\usepackage{amssymb,hyperref}

\theoremstyle{plain}
\newtheorem{thm}{Theorem}
\newtheorem{lem}{Lemma}
\newtheorem{cor}{Corollary}

\theoremstyle{definition}
\newtheorem{ex}{Example}

\renewcommand{\Re}{\mathrm{Re}}

\title{New extensions for a theorem by Mocanu}

\author{Hitoshi Shiraishi}
\address{Hitoshi Shiraishi \newline
Department of Mathematics \newline
Kinki University \newline
Higashi-Osaka, Osaka 577-8502, Japan}
\email{shiraishi@math.kindai.ac.jp}

\author{Shigeyoshi Owa}
\address{Shigeyoshi Owa \newline
Department of Mathematics \newline
Kinki University \newline
Higashi-Osaka, Osaka 577-8502, Japan}
\email{owa@math.kindai.ac.jp}


\date{}

\begin{document}

\begin{abstract}
For analytic functions $f(z)$ in the open unit disk $\mathbb{U}$ with $f(0)=f'(0)-1=f''(0)=0$,
P. T. Mocanu (Mathematica (Cluj), {\bf 42}(2000)) has considered some sufficient arguments of $f'(z)+zf''(z)$ for $|\arg(zf'(z)/f(z))|<\pi\mu/2$.
The object of the present paper is to discuss those probrems for $f(z)$ with $f''(0)=f'''(0)=\ldots=f^{(n)}(0)=0$ and $f^{(n+1)}(0) \ne 0$.
\end{abstract}

\maketitle

\section{Introduction}

\

Let $\mathcal{A}_n$ denote the class of functions
$$
f(z)=z+a_{n+1}z^{n+1}+a_{n+2}z^{n+2}+ \ldots
\qquad(n=1,2,3,\ldots)
$$
that are analytic in the open unit disk $\mathbb{U}=\{z \in \mathbb{C}:|z|<1\}$
and $\mathcal{A}=\mathcal{A}_1$.

Also,
let $\mathcal{H}[1,n]$ denote the class of functions $f(z)$ of the form
$$
f(z)
= 1 + \sum_{k=n}^{\infty} a_k z^k
\qquad (n=1,2,3,\ldots)
$$
which are analytic in $\mathbb{U}$.

Further,
let the class $\mathcal{STS}(\mu)$ of $f(z) \in \mathcal{A}_n$ be defined by
$$
\mathcal{STS}(\mu)=\left\{f(z)\in\mathcal{A}_n:\left|\arg\left(\frac{zf'(z)}{f(z)}\right)\right|<\frac{\pi}{2}\mu,\ 0<\mu\leqq1\right\}
$$
and $\mathcal{S}^*=\mathcal{STS}(1)$. This class $\mathcal{STS}(\mu)$ was considered by Shiraishi and Owa \cite{m1ref4}.

Let $f(z)$ and $g(z)$ be analytic in $\mathbb{U}$.
Then $f(z)$ is said to be subordinate to $g(z)$ if there exists an analytic function $w(z)$ in $\mathbb{U}$ satisfying $w(0)=0$, $|w(z)| < 1\,\,(z\in\mathbb{U})$ and $f(z)=g(w(z))$.
We denote this subordination by
$$
f(z)
\prec g(z)
\qquad (z\in\mathbb{U}).
$$

\

\section{Lemmas}

\

We need the following lemmas to consider our main results.

\

\begin{lem} \label{d6lem1} \quad
Let $n$ be a positive integer,
$\lambda>0$,
and let $\beta_0 = \beta_0(\lambda,n)$ be the root of the equation
$$
\beta\pi
= \frac{3}{2} \pi - \arctan( n \lambda \beta).
$$

In addition,
let 
$$
\alpha
= \alpha(\beta,\lambda,n)
= \beta + \frac{2}{\pi}\arctan(n\lambda\beta)
$$
for $0<\beta\leqq\beta_0$.

If $p(z) \in \mathcal{H}[1,n]$ and
$$
p(z) + \lambda zp'(z)
\prec \left( \frac{1+z}{1-z} \right)^\alpha
\qquad (z\in\mathbb{U}),
$$
then
$$
p(z)
\prec \left( \frac{1+z}{1-z} \right)^\beta
\qquad (z\in\mathbb{U}).
$$
\end{lem}

\

\begin{lem} \label{d6lem2} \quad
Let $q(z)$ be the convex function in $\mathbb{U}$,
with $q(0)=1$ and $\Re(q(z))>0$ for $\mathbb{U}$.
Let the function $h(z)$ be given by
$$
h(z)
= (q(z))^2+nzq'(z)
\qquad (z\in\mathbb{U}).
$$

If $p(z) \in \mathcal{H}[1,n]$ and
$$
(p(z))^2+zp'(z)
\prec h(z)
\qquad (z\in\mathbb{U}),
$$
then $p(z) \prec q(z)$ and this is sharp.
\end{lem}

\

The above lemmas were given by Mocanu \cite{d6ref3}.

Applying Lemma \ref{d6lem2},
we obtain the following lemma.

\

\begin{lem} \label{d6lem3} \quad
If $p(z) \in \mathcal{H}[1,n]$ satisfies
$$
|\arg((p(z))^2+zp'(z))|
<\phi(\mu)
\qquad (z\in\mathbb{U})
$$
for
$$
\phi(\mu)
= \frac{\pi}{2}(\mu+1) - \arctan \frac{\cos \frac{\mu\pi}{2}}{\sin \frac{\mu\pi}{2} + \frac{n\mu}{1-\mu} \left( \frac{1-\mu}{1+\mu} \right)^\frac{1+\mu}{2} }
$$
and $0<\mu\leqq1$,
then
$$
|\arg (p(z))|
< \frac{\pi}{2} \mu
\qquad (z\in\mathbb{U}).
$$
\end{lem}

\

\begin{proof} \qquad
Let us define the function $q(z)$ by
\begin{equation} \label{d6lem3eq1}
q(z)
= \left( \frac{1+z}{1-z} \right)^\mu
\qquad (z\in\mathbb{U})
\end{equation}
for $0<\mu\leqq1$ and the function $h(z)$ by 
$$
h(z)
= (q(z))^2+nzq'(z)
\qquad (z\in\mathbb{U}).
$$

Then the function $q(z)$ is convex in $\mathbb{U}$
with $\Re(q(z))>0$,
$h(z)$ is univalent in $\mathbb{U}$
and $h(\mathbb{U})$ is the symmetric domain
with respect to the real axis.

If we set $z=\exp(i\theta)$,
$0\leqq\theta<\pi$
and $x=\cot\dfrac{\theta}{2}$,
then $x\geqq0$,
$z=\dfrac{ix-1}{ix+1}$
and $q(z)=(ix)^\mu$.
Hence
$$
h(e^{i\theta})
= (ix)^{\mu-1}H(x),
$$
where
$$
H(x)
= (ix)^{\mu+1}-\frac{n}{2}\mu(1+x^2).
$$

Noting that $\cos\dfrac{(\mu+1)\pi}{2}=-\sin\dfrac{\mu\pi}{2}$
and $\sin\dfrac{(\mu+1)\pi}{2}=\cos\dfrac{\mu\pi}{2}$,
we see that
$$
H(x)
= P(x)+iQ(x),
$$
where
$$
\left\{
\begin{array}{l}
P(x)
= -\sin\dfrac{\mu\pi}{2} x^{\mu+1}-\dfrac{n}{2}\mu(1+x^2) \\
Q(x)
= \cos\dfrac{\mu\pi}{2} x^{\mu+1}.
\end{array}
\right.
$$

Let
$$
\varphi(\mu)
= \min \{ \arg (H(x)) : x\geqq0 \}
$$
and
\begin{equation} \label{d6lem3eq2}
\phi(\mu)
= \varphi(\mu) + \frac{\pi}{2}(\mu-1).
\end{equation}

From (\ref{d6lem3eq1})
we deduce
$$
\arg (h(e^{i\theta}))
\geqq \phi(\mu).
$$

Since
$$
G(x)
= Q'(x)P(x)-P'(x)Q(x)
= \frac{n}{2} \mu \cos\dfrac{\mu\pi}{2} x^\mu ((1-\mu)x^2-(1+\mu))
= 0
$$
has the root $x_0=\left( \dfrac{1-\mu}{1+\mu} \right)^\frac{1}{2}$ and
$$
\frac{Q(x_0)}{P(x_0)}
= - \frac{\cos \frac{\mu\pi}{2}}{\sin \frac{\mu\pi}{2} + \frac{n\mu}{1-\mu} \left( \frac{1-\mu}{1+\mu} \right)^\frac{1+\mu}{2} },
$$
we deduce
\begin{equation} \label{d6lem3eq3}
\varphi(\mu)
= \pi - \arctan \frac{\cos \frac{\mu\pi}{2}}{\sin \frac{\mu\pi}{2} + \frac{n\mu}{1-\mu} \left( \frac{1-\mu}{1+\mu} \right)^\frac{1+\mu}{2} }.
\end{equation}

Hence
$$
|\arg (h(e^{i\theta}))|
\geqq \phi(\mu)
\qquad (-\pi<\theta<\pi)
$$
where $\phi(\mu)$ is given by (\ref{d6lem3eq2}) and (\ref{d6lem3eq3}).

From the assumption,
$$
(p(z))^2+zp'(z)
\prec h(z)
\qquad (z\in\mathbb{U}).
$$

Hence by Lemma \ref{d6lem2}
we deduce $p(z) \prec q(z)$.

So,
we obtain
$$
|\arg (p(z))|
< \frac{\pi}{2} \mu
\qquad (z\in\mathbb{U}).
$$
\end{proof}

\

\section{Main results}

\

Using Lemma \ref{d6lem1} and Lemma \ref{d6lem3},
we get the following result.

\

\begin{thm} \label{d6thm1} \quad
If $f(z) \in \mathcal{A}_n$ satisfies
$$
|\arg(f'(z)+zf''(z))|
< \frac{\pi}{2}\alpha
\qquad (z\in\mathbb{U})
$$
for
$$
\alpha
= \beta+\frac{2}{\pi} \arctan (n\beta),
$$
$$
\beta
= \gamma+\frac{2}{\pi} \arctan (n\gamma),
$$
$$
\frac{\pi}{2}(\alpha+\gamma)
\leqq \phi(\mu)
$$
and
$$
\phi(\mu)
= \frac{\pi}{2}(\mu+1) - \arctan \frac{\cos \frac{\mu\pi}{2}}{\sin \frac{\mu\pi}{2} + \frac{n\mu}{1-\mu} \left( \frac{1-\mu}{1+\mu} \right)^\frac{1+\mu}{2} }
$$
with some real $\alpha,\gamma>0$ and $0<\mu\leqq1$,
then $f(z)\in\mathcal{STS}(\mu)$.
\end{thm}

\

\begin{proof} \qquad
By using Lemma \ref{d6lem1} with $\lambda=1$
we deduce
$$
|\arg (f'(z))|
< \frac{\pi}{2}\beta
\qquad (z\in\mathbb{U}).
$$
with
$$
\alpha
= \beta+\frac{2}{\pi} \arctan (n\beta).
$$

Using again Lemma \ref{d6lem1},
we get
$$
\left| \arg \left(\frac{f(z)}{z}\right) \right|
< \frac{\pi}{2}\gamma
\qquad (z\in\mathbb{U}),
$$
where $\gamma$ is the solution of the equation
$$
\beta
= \gamma+\frac{2}{\pi} \arctan (n\gamma).
$$

If we set $p(z)=\dfrac{zf'(z)}{f(z)}$ and $P(z)=\dfrac{f(z)}{z}$,
then we have $p(z) \in \mathcal{H}[1,n]$ and
$$
f'(z)+zf''(z)
= P(z)((p(z))^2+zp'(z))
\qquad (z\in\mathbb{U}),
$$
where
$$
|\arg (P(z))|
< \frac{\pi}{2}\gamma
\qquad (z\in\mathbb{U}).
$$

It follows that
$$
|\arg((p(z))^2+zp'(z))|
\leqq |\arg(f'(z)+zf''(z))| + |\arg (P(z))|
< \frac{\pi}{2}(\alpha+\gamma).
$$

For the condition of $\phi(\mu)$,
we deduce that
$$
|\arg((p(z))^2+zp'(z))|
< \phi(\mu)
\qquad (z\in\mathbb{U})
$$
implies by means of Lemma \ref{d6lem3}.
$$
|\arg (p(z))|
< \dfrac{\pi}{2}\mu
\qquad (z\in\mathbb{U}).
$$
\end{proof}

\

We consider an example for Theorem \ref{d6thm1}.

\

\begin{ex} \label{d6ex1} \quad
Let us consider the function
$$
f(z)
= z + \sin \frac{\pi\alpha}{2} z^{n+1}
\qquad(z\in\mathbb{U})
$$
with $0<\alpha\leqq1$.
If we put
$$
\mu
= \frac{2}{\pi} \arcsin \frac{n(n+1)\sin\frac{\pi\alpha}{2}}{(n+1)^3-\sin^2\frac{\pi\alpha}{2}},
$$
the function $f(z)$ satisfies the condition of Theorem \ref{d6thm1}.

Because differenciating the function $f(z)$,
we obtain
\begin{align*}
\frac{zf'(z)}{f(z)}
&= \frac{(n+1)(n+1+\sin\frac{\pi\alpha}{2}z^n)}{(n+1)^2+\sin\frac{\pi\alpha}{2}z^n} \\
&= n+1-\frac{n}{1+\frac{\sin\frac{\pi\alpha}{2}}{(n+1)^2}z^n}
\qquad(z\in\mathbb{U})
\end{align*}
and therefore,
$$
\left| \arg \left(\frac{zf'(z)}{f(z)}\right) \right|
< \arcsin \frac{n(n+1)\sin\frac{\pi\alpha}{2}}{(n+1)^3-\sin^2\frac{\pi\alpha}{2}}
\qquad(z\in\mathbb{U}).
$$
\end{ex}

\

If we fix one of the values for $\alpha$, $\beta$ or $\gamma$ in Theorem \ref{d6thm1},
then we can obtain others.
For example,
if we put $n=2$, $\alpha=1$ and $\mu=\dfrac{1}{2}$,
then we get $\beta=\dfrac{1}{2}$, $\gamma=0.227\ldots$ and the following result due to Mocanu \cite{d6ref3}.

\

\begin{cor} \label{d6cor1} \quad
If $f(z) \in \mathcal{A}_2$ satisfies
$$
\Re(f'(z)+zf''(z))
> 0
\qquad (z\in\mathbb{U}),
$$
then
$$
\left| \arg \left(\frac{zf'(z)}{f(z)}\right) \right|
< \frac{\pi}{4}
\qquad (z\in\mathbb{U}).
$$
\end{cor}

\

Moreover,
putting $n=2$, $\alpha=\dfrac{3}{2}$ and $\mu=1$,
we have Corollary \ref{d6cor2} given by Mocanu \cite{d6ref3}.

\

\begin{cor} \label{d6cor2} \quad
If $f(z) \in \mathcal{A}_2$ satisfies
$$
|\arg(f'(z)+zf''(z))|
< \frac{3}{4}\pi
\qquad (z\in\mathbb{U}),
$$
then
$$
\Re\left(\frac{zf'(z)}{f(z)}\right)
>0
\qquad (z\in\mathbb{U}).
$$
\end{cor}

\

Futhermore,
putting $n=1$,
we derive Corollary \ref{d6cor3} and Corollary \ref{d6cor4} which were showed by Mocanu \cite{d6ref2}.

\

\begin{cor} \label{d6cor3} \quad
If $f(z) \in \mathcal{A}$ satisfies
$$
\Re(f'(z)+zf''(z))
> 0
\qquad (z\in\mathbb{U}),
$$
then
$$
\left| \arg \left(\frac{zf'(z)}{f(z)}\right) \right|
< \frac{\pi}{3}
\qquad (z\in\mathbb{U}).
$$
\end{cor}

\

\begin{cor} \label{d6cor4} \quad
If $f(z) \in \mathcal{A}$ satisfies
$$
|\arg(f'(z)+zf''(z))|
< \frac{2}{3}\pi
\qquad (z\in\mathbb{U}),
$$
then
$$
\Re\left(\frac{zf'(z)}{f(z)}\right)
>0
\qquad (z\in\mathbb{U}).
$$
\end{cor}

\

\section{Integral version of the results}

\

Let us define the function
$$
F(z)
= \int^z_0 \frac{f(t)}{t} dt
\qquad (z\in\mathbb{U})
$$
for $f(z)\in\mathcal{A}_n$.
This integral operator $F(z)$ is given by Alexander \cite{d6ref1} and is said to be Alexander integral operator.

For this Alexander integral operator for $f(z)$,
we derive

\

\begin{thm} \label{d6thm2} \quad
If $f(z) \in \mathcal{A}_n$ satisfies
$$
|\arg (f'(z))|
< \frac{\pi}{2}\alpha
\qquad (z\in\mathbb{U})
$$
for
$$
\alpha
= \beta+\frac{2}{\pi} \arctan (n\beta),
$$
$$
\beta
= \gamma+\frac{2}{\pi} \arctan (n\gamma),
$$
$$
\frac{\pi}{2}(\alpha+\gamma)
\leqq \phi(\mu)
$$
and
$$
\phi(\mu)
= \frac{\pi}{2}(\mu+1) - \arctan \frac{\cos \frac{\mu\pi}{2}}{\sin \frac{\mu\pi}{2} + \frac{n\mu}{1-\mu} \left( \frac{1-\mu}{1+\mu} \right)^\frac{1+\mu}{2} }
$$
with some real $\alpha,\gamma>0$ and $0<\mu\leqq1$,
then the Alexander integral operator $F(z)$ of $f(z)$ belongs to the class $\mathcal{STS}(\mu)$.
\end{thm}

\

The proof of Theorem \ref{d6thm2} follows by replacing $f(z)$ with $F(z)$ in Theorem \ref{d6thm1}.

\

\end{document}